# A COMPARISON OF EIGENVALUE-BASED ALGORITHMS AND THE GENERALIZED LANCZOS TRUST-REGION ALGORITHM FOR SOLVING THE TRUST-REGION SUBPROBLEM[*]

ZHONGXIAO JIA[†] AND FA WANG[‡]

**Abstract.** Solving the trust-region subproblem (TRS) plays a key role in numerical optimization and many other applications. Based on a fundamental result that the solution of TRS of size $n$ is mathematically equivalent to finding the rightmost eigenpair of a certain matrix pair of size $2n$, eigenvalue-based methods are promising due to their simplicity. For $n$ large, the implicitly restarted Arnoldi (IRA) and refined Arnoldi (IRRA) algorithms are well suited for this eigenproblem. For a reasonable comparison of overall efficiency of the algorithms for solving TRS directly and eigenvalue-based algorithms, a vital premise is that the two kinds of algorithms must compute the approximate solutions of TRS with (almost) the same accuracy, but such premise has been ignored in the literature. To this end, we establish close relationships between the two kinds of residual norms, so that, given a stopping tolerance for IRA and IRRA, we are able to determine a reliable one that GLTR should use so as to ensure that GLTR and IRA, IRRA deliver the converged approximate solutions with similar accuracy. We also make a convergence analysis on the residual norms by the Generalized Lanczos Trust-Region (GLTR) algorithm for solving TRS directly, the Arnoldi method and the refined Arnoldi method for the equivalent eigenproblem. A number of numerical experiments are reported to illustrate that IRA and IRRA are competitive with GLTR and IRRA outperforms IRA.

**Key words.** trust-region subproblem, eigenvalue problem, implicitly restarted Arnoldi algorithm, implicitly restarted refined Arnoldi algorithm, GLTR algorithm, Krylov subspace, residual norm, stopping criterion

**AMS subject classifications.** 90C20, 90C30, 65K05, 49M37

**1. Introduction.** The trust-region subproblem (TRS) is

$$(1.1) \qquad \min_{\|s\|_B \leq \Delta} q(s) \quad \text{with} \quad q(s) = g^T s + \frac{1}{2} s^T A s,$$

where $A \in \mathbb{R}^{n \times n}$ is symmetric, $B \in \mathbb{R}^{n \times n}$ is symmetric positive definite, the nonzero vector $g \in \mathbb{R}^n$, $\Delta > 0$ is the trust-region radius, and the norm $\|\cdot\|_B$ is the $B$-norm defined by $\|s\|_B = \sqrt{s^T B s}$. The matrix $B$ is not necessarily equal to the identity $I$, and it is often constructed to impose a smoothness condition on a solution to (1.1) for the ill-posed problem and to incorporate scaling of variables in optimization [22]. For instance, it is argued in [2] that a good choice is $B = J^{-T} J^{-1}$ for some invertible matrix $J$ or the Hermitian polar factor of $A$ [13].

Problem (1.1) arises from many applications, e.g., from nonlinear numerical optimization [2, 21], where $q(s)$ is a quadratic model of $\min f(s)$ at the current approximate solution, $A$ is Hessian and $g$ is the gradient of $f$ at the current approximate solution, and many others, e.g., Tikhonov regularization of ill-posed problems [22, 23], graph partitioning problems [12], the constrained eigenvalue problem [6], and the Levenberg–Marquardt algorithm for solving nonlinear least squares problems [21].

The following theorem [2, 20] provides a theoretical basis for many TRS solution algorithms and give necessary and sufficient conditions, called the optimal conditions,

---

[*]This work was supported in part by the National Natural Science Foundation of China (No. 11771249)

[†]Corresponding author. Department of Mathematical Sciences, Tsinghua University, 100084 Beijing, China. (jiazx@tsinghua.edu.cn)

[‡]Department of Mathematical Sciences, Tsinghua University, 100084 Beijing, China. (wangfa15@mails.tsinghua.edu.cn)

for the solution of TRS (1.1).

THEOREM 1.1. *A vector $s_{opt}$ is a solution to* (1.1) *if and only if there exists the optimal Lagrangian multiplier $\lambda_{opt} \geq 0$ such that*

$$\|s_{opt}\|_B \leqslant \Delta, \tag{1.2}$$
$$(A + \lambda_{opt} B)s_{opt} = -g, \tag{1.3}$$
$$\lambda_{opt}(\Delta - \|s_{opt}\|_B) = 0, \tag{1.4}$$
$$A + \lambda_{opt} B \succeq 0, \tag{1.5}$$

*where the notation $\succeq 0$ indicates that a symmetric matrix is semi-positive definite.*

TRS algorithms have been extensively studied for a few decades, and a number of algorithms have been available. On the basis of [1], the authors in [19] have classified the known algorithms into four categories: accurate methods for small to medium sized problems [20], accurate methods for large sparse problems [4, 10, 11, 24, 22, 23, 27], approximate methods [28, 31] which are designed to approximately solve (1.1), and eigenvalue-based methods [1, 6], which are suitable for both small to medium sized problems and large sparse problems.

Over the past two decades, the GLTR method [8] has become a commonly used algorithm for large scale TRS's. It consists of two phases: In the first phase, it starts with the Truncated Conjugate Gradient (TCG) algorithm with the zero initial guess [28, 31]. When $A$ is positive definite and $\|A^{-1}g\|_B \leq \Delta$, which corresponds to $\lambda_{opt} = 0$, the algorithm ultimately returns a converged approximate solution to $s_{opt} = -A^{-1}g$ and stops. If the TCG iterate exceeds the trust-region boundary or a negative curvature is encountered, which means that the trust-region constraint is activated and corresponds to a semi-positive definite or indefinite $A$ but $\lambda_{opt} > 0$, the GLTR method enters the second phase and switches to the Lanczos method that accurately solves a sequence of the projected problems restricted to some Krylov subspaces, and it proceeds in such a way until the approximate solution converges to $s_{opt}$. The GLTR method may require many iterations, so that storing the basis vectors may become unacceptable. One may store them in secondary memory and calls them when forming GLTR iterates, but this way is time consuming. Besides, numerical difficulties such as loss of orthogonality among the Lanczos basis vectors make the algorithm intricate. To this end, Zhang et al. [33, 34] have recently proposed a complicated restarted GLTR algorithm, which needs further explorations and developments.

A potentially great breakthrough has come to the solution of TRS's when Adachi et al. [1] have proved that TRS (1.1) can be treated by solving the generalized eigenvalue problem of a certain matrix pair of order $2n$ if the solution $s_{opt}$ to the TRS (1.1) is in the easy case and reaches the trust-region boundary. Precisely, they have shown that a solution to TRS (1.1) can be determined by the rightmost eigenvalue and the associated eigenvector of some $2n \times 2n$ matrix pair. Therefore, computing the solution $s_{opt}$ of TRS (1.1) amounts to computing the rightmost eigenpair of a generalized eigenvalue problem, for which a number of standard numerical methods are applicable [7, 29], depending on the size of $n$. Concretely, for $A$ large, Adachi et al. have used the implicitly restarted Arnoldi (IRA) algorithm proposed by Sorensen [26], i.e., the Matlab function eigs, to solve the generalized eigenvalue problem. IRA is simple to implement, and has the advantage of low storage and ingenious restarting strategy. Numerical experiments in [1, 34] have indicated that IRA can be competitive with the GLTR algorithm in terms of CPU time and accuracy and may thus be a promising algorithm for large scale TRS's. However, the comparison standards used there have serious oversights or gaps. They set stopping tolerances at will for residual norms by the algorithms applied to (1.1) directly, such as GLTR, and the ones for the eigenvalue-based algorithms, such as IRA. Such a way causes unreasonable or incorrect efficiency comparisons because the same or similar stopping criteria for the two kinds of different algorithms that solve (1.1) and its mathematically equivalent eigenproblem do not mean that they deliver converged approximations of $s_{opt}$ with similar accuracy. Therefore, we should find relationships between them, relate them and set them correctly in order to make two kinds of algorithms give almost the same solution accuracy for TRS (1.1) when making an efficiency comparison.

For large scale matrix eigenproblems, it has been known from, e.g., [14, 18, 29], that Arnoldi algorithms, such as IRA, may have convergence problems [14, 18] for computing eigenvectors. Precisely, the approximate eigenvectors, i.e., the Ritz vectors, may converge erratically and even may not converge even though the corresponding Ritz values converge. Jia [15, 16] has corrected this problem by proposing the refined Rayleigh–Ritz method that uses new approximate eigenvectors, called refined Ritz vectors or refined eigenvector approximations, to approximate the desired eigenvectors; see [16, 29] for details. More generally, given a general projection subspace, the refined Ritz vector is more accurate than the Ritz vector [17]. Particularly, for the refined Arnoldi method [15], Jia [16] has developed an implicitly restarted refined Arnoldi (IRRA) algorithm, which applies the same implicitly restarting scheme as in IRA to the refined Arnoldi method but uses the new refined shifts to construct a better projection subspace at each restart than IRA does, so that IRRA may outperform IRA substantially.

In this paper, we extend the refined Arnoldi method [15] and IRRA [18] for the standard eigenvalue problems to generalized ones. Then our main goal is to make extensive experiments on TRS's with $B = I$ and $B \neq I$ and to show whether or not IRA and IRRA are practically competitive with the GLTR algorithm without restart. For this purpose, we need to settle down two important issues. One of them is to select a fair comparison measure. As is well known, in exact arithmetic, the Lanczos process in GLTR exploits elegant three term recurrences to compute an orthonormal or $B$-orthonormal basis of the underlying Krylov subspace with $B = I$ or $B \neq I$, but the Arnoldi process in IRA and IRRA uses long recurrences to do similar things. For subspaces of small dimension, due to algorithmic features, the number of matrix-vector products (MVs) dominates the efficiency of the GLTR method, the Arnoldi method and the refined Arnoldi method for subspaces of small dimension. As a result, the total number of MVs for GLTR, IRA and IRRA is a reasonable measure to judge their overall efficiency when subspace dimension for IRA and IRRA is small. There are certainly other measures, such as the CPU time. But the reliability of CPU time is much involved, and it strongly depends on programming optimization, programming languages and computing environment, etc. In any event, we must stress that MVs is only a rough measure of overall performance, and we ignore complicated effects of several other factors, e.g., the sparsity of $A$ and $B$, reorthogonalization in the Lanczos process and the Arnoldi process, and calls of the Lanczos basis vectors from secondary storages when iterations used by GLTR are too many and secondary storages have to be exploited.

The other important issue that deserves enough attention is to establish correct and fair stopping criteria when comparing the overall efficiency of GLTR and IRA, IRRA. This issue has received no attention in the literature, and one commonly terminates them by simply setting user-prescribed tolerances for respective residual norms in an irrelevant way [1, 34]. This will cause unreasonable or incorrect comparisons because the residual norms themselves for GLTR and IRA, IRRA are completely different and are generally not of similar or comparable size. That is, the same tolerances for the two kinds of algorithms do not mean that the TRS residual norm of approximate solution recovered from the eigenvalue-based algorithms is the same as or very comparable to the residual norm by algorithms for solving TRS directly, causing that two kinds of algorithms do not compute the converged approximate solutions of TRS with very similar accuracy.

In view of the above, as the premise of efficiency comparison of two kinds of algorithms, we should fix one stopping criterion for one kind of algorithm as the reference standard and then determine how small the other kind should be, so that we can finally compute approximate solutions with almost the same accuracy by the two kinds of algorithms. In this paper, we shall establish a concise relationship between the different stopping criteria for GLTR and IRA, IRRA. Based on it, we prescribe the stopping criteria for IRA and IRRA, then decide the one for GLTR, so that we can make a reasonable efficiency comparison. In the meantime, we will theoretically investigate how fast GLTR, the Arnoldi method and the refined Arnoldi method converge by analyzing the residual norms obtained by GLTR and by IRA and IRRA for the equivalent matrix eigenproblem.

This paper is organized as follows. In Section 2, we give some preliminaries and review the GLTR method. In Section 3, we introduce the equivalence of the solutions of (1.1) and the generalized eigenvalue problem of a certain $2n \times 2n$ matrix pair. We describe a general framework of an eigenvalue-based algorithm in [1], and introduce IRA [26]. We then propose the refined Arnoldi method and develop the corresponding IRRA algorithm [16] in the context of generalized eigenvalue problems. Section 4 is devoted to deriving relationships between the TRS residual norms and residual norms of the equivalent generalized eigenvalue problem. In Section 5, we investigate the convergence of residual norms obtained by GLTR, the Arnoldi method and the refined Arnoldi method for the equivalent eigenproblem of TRS. In Section 6, we report numerical experiments to compare IRA and IRRA with GLTR, showing that they are competitive with GLTR. Meanwhile, we compare IRRA with IRA and demonstrate that IRRA is more efficient than IRA. Finally, we conclude the paper in Section 7.

Throughout this paper, denote by the superscript $T$ the transpose of a matrix or vector, by $\|\cdot\|$ the 2-norm of a matrix or vector, by $I$ the identity matrix with order clear from the context, and by $e_i$ the $i$th column of $I$. All vectors are column vectors and are typeset in lowercase letters.

## 2. Preliminaraies.

**2.1. A solution to** (1.1). As we can see from Theorem 1.1, if (1.1) has no solution with $\|s_{opt}\|_B = \Delta$, then $A$ is positive definite and $s_{opt} = -A^{-1}g$ with $\|s_{opt}\|_B < \Delta$ and $\lambda_{opt} = 0$. Otherwise, $\|s_{opt}\|_B = \Delta$ and $\lambda_{opt} > 0$, the solution $s_{opt}$ to TRS (1.1) is unique and $s_{opt} = -(A + \lambda_{opt}B)^{-1}g$. All these correspond to the so-called "easy case" [2, 8, 20, 21] or "nondegenerate case" [2, 8, 20, 21].

If $A$ is semi-positive definite or indefinite and

$$g \perp \mathcal{N}(A - \alpha_n B),$$

the null space of $A - \alpha_n B$, where $\alpha_n$ is the leftmost eigenvalue of the matrix pair $(A, B)$, then we have the following definition [2, 8, 21].

DEFINITION 2.1 (Hard Case). *The solution of TRS* (1.1) *is a hard case if $g$ is*

*orthogonal to the eigenspace corresponding to the eigenvalue $\alpha_n$ of $A - \lambda B$, and the optimal Lagrangian multiplier is $\lambda_{opt} = -\alpha_n$.*

In the hard case, TRS (1.1) may have multiple optimal solutions [21, p.87-88], which must reach the trust-region boundary and can be expressed by

$$(2.1) \qquad s_{opt} = -(A - \alpha_n B)^\dagger g + \eta u_n,$$

where $u_n \in \mathcal{N}(A - \alpha_n B)$, $\|(A - \alpha_n B)^\dagger g\|_B < \Delta$, and the superscript $\dagger$ denotes the Moore-Penrose generalized inverse.

It is known that the scalar $\eta$ satisfies

$$(2.2) \qquad \eta^2 = \Delta^2 - \|(A - \alpha_n B)^\dagger g\|_B^2 > 0$$

when $\|u_n\|_B = 1$ because of the $B$-orthogonality of $(A - \alpha_n B)^\dagger g$ and $u_n$.

From (2.1), we see that there are infinitely many $s_{opt}$ whenever the eigenvalue $\alpha_n$ of the matrix pair $(A, B)$ is multiple, since $\mathcal{N}(A - \alpha_n B)$ is bigger than one. We not only need to find the minimum $B$-norm solution to a symmetric semi-positive definite singular system but also need to compute the eigenspace of the matrix pair $(A, B)$ associated with the eigenvalue $\alpha_n$. The hard case must be treated separately and has been studied for years; see, e.g., [1, 5, 8, 20, 21, 24].

As has been addressed in [2], the hard case rarely occurs in practice, as it requires that both $A$ be semi-definite or indefinite and $g$ be orthogonal to $\mathcal{N}(A - \alpha_n B)$. In the sequel, we are only concerned with the easy case.

**2.2. The generalized Lanczos trust-region (GLTR) method [8].** Without loss of generality, we always assume that $\|s_{opt}\|_B = \Delta$ throughout the paper for otherwise $A$ must be symmetric positive definite and the Conjugate Gradient (CG) method can solve (1.1). In this subsection, we briefly review the GLTR method [8].

At iteration $k$, the GLTR method exploits the symmetric Lanczos process to generate a $B$-orthonormal basis $\{p_i\}_{i=0}^k$ of the underlying $(k+1)$-dimensional Krylov subspace

$$\mathcal{K}_k(B^{-1}g, B^{-1}A) = span\{B^{-1}g, (B^{-1}A)B^{-1}g, \ldots, (B^{-1}A)^k B^{-1}g\},$$

which can be written in matrix form

$$(2.3) \qquad AP_k = BP_k T_k + \beta_{k+1} B p_{k+1} e_{k+1}^T,$$

$$(2.4) \qquad P_k^T g = \beta_0 e_1, \ \beta_0 = \|B^{-1}g\|_B = \|B^{-\frac{1}{2}}g\|,$$

$$(2.5) \qquad g = \beta_0 B q_0,$$

where $P_k = (p_0, p_1, \ldots, p_k)$ is $B$-orthonormal, i.e., $P_k^T B P_k = I$, the matrix

$$(2.6) \qquad T_k = P_k^T A P_k = \begin{pmatrix} \delta_0 & \beta_1 & & & \\ \beta_1 & \delta_1 & \ddots & & \\ & \ddots & \ddots & \ddots & \\ & & \ddots & \delta_{k-1} & \beta_k \\ & & & \beta_k & \delta_k \end{pmatrix} \in \mathbb{R}^{(k+1) \times (k+1)}$$

is symmetric tridiagonal, and $e_{k+1}$ is the last column of the identity matrix of order $k+1$. The GLTR method projects TRS (1.1) onto the $(k+1)$-dimensional Krylov subspace $\mathcal{K}_k(B^{-1}g, B^{-1}A)$ and yields the reduced TRS:

$$\min_{\|h\|\leq \Delta} \phi(h) = \beta_0 e_1^T h + \frac{1}{2} h^T T_k h. \tag{2.7}$$

From Theorem 1.1, the vector $h^{(k)}$ is a solution to (2.7) if and only if there exists the optimal Lagrangian multiplier $\lambda^{(k)} \geq 0$ such that

$$\|h^{(k)}\| \leq \Delta,$$
$$(T_k + \lambda^{(k)} I) h^{(k)} = -\beta_0 e_1,$$
$$\lambda^{(k)}(\Delta - \|h^{(k)}\|) = 0,$$
$$T_k + \lambda^{(k)} I \succeq 0.$$

GLTR computes the iterate $s_{\text{gltr}}^{(k)} = P_k h_k$ to approximate $s_{opt}$, which solves the projected TRS

$$\min_{s \in \mathcal{K}_k(B^{-1}g, B^{-1}A), \|s\|_B \leq \Delta} q(s). \tag{2.8}$$

The GLTR method exploits the Moré-Sorensen method [10] to solve the reduced TRS (2.7). It is shown in [8] that the $B^{-1}$-residual norm of the approximate solution $(\lambda^{(k)}, s_{\text{gltr}}^{(k)})$ of (1.3) satisfies

$$\|(A + \lambda^{(k)} B) s_{\text{gltr}}^{(k)}\|_{B^{-1}} = \beta_{k+1} |e_{k+1}^T h^{(k)}|, \tag{2.9}$$

which is used as an efficient stopping criterion without forming $s_{\text{gltr}}^{(k)}$ explicitly until the residual norm drops below a user-prescribed tolerance.

Importantly, if $\|s_{opt}\|_B = \Delta$, we must have $\|s_{\text{gltr}}^{(k)}\|_B = \Delta$ as $k$ increases [8, 33]. In fact, it has been shown in [19] that $\|s_{\text{gltr}}^{(k)}\|_B = \Delta$ after very few iterations. This has been numerically confirmed to occur for $k = 1$ in [19, 33].

In order to get an accurate approximation, one needs to expand the Krylov subspace successively and store $P_k$ so as to form the iterate $s_{\text{gltr}}^{(k)}$. In finite precision arithmetic, numerical difficulties such as loss of orthogonality and excessive storage may limit the applicability of the method. A complicated restarted variant of GLTR has been proposed in [34] to reduce storage, and it has been numerically justified to be comparable to GLTR in terms of MVs on some problems.

**3. The equivalence of TRS and a certain matrix eigenvalue problem and the IRA, IRRA algorithms.**

**3.1. The equivalence of TRS (1.1) and a matrix eigenvalue problem.** Define

$$M = \begin{pmatrix} -A & \frac{gg^T}{\Delta^2} \\ B & -A \end{pmatrix} \in \mathbb{R}^{2n \times 2n}, \quad \widetilde{B} = \begin{pmatrix} B & 0 \\ 0 & B \end{pmatrix} \in \mathbb{R}^{2n \times 2n}. \tag{3.1}$$

Adachi et al. [1] prove that if $\|s_{opt}\|_B = \Delta$ then the solution of TRS (1.1) can be equivalently formulated as the rightmost eigenpair of the matrix pair $(M, \widetilde{B})$, as shown below.

THEOREM 3.1 ([1]). *Let $(\lambda_{opt}, s_{opt})$ satisfy Theorem 1.1 with $\|s_{opt}\|_B = \Delta$. Then the rightmost eigenvalue $\mu_1$ of $(M, \widetilde{B})$ is real and simple, and the optimal Lagrangian* multiplier $\lambda_{opt} = \mu_1$. Let $y = (y_1^T, y_2^T)^T$ be an eigenvector of $(M, \widetilde{B})$ associated with $\mu_1$, and suppose that $g^T y_2 \neq 0$. Then the unique TRS solution is

$$s_{opt} = -\frac{\Delta^2}{g^T y_2} y_1. \tag{3.2}$$

REMARK 3.1. *In finite precision arithmetic, it is shown in [1] that $s_{opt}$ can be computed stably by the formula*

$$s_{opt} = -sign(g^T y_2) \Delta \frac{y_1}{\|y_1\|_B}. \tag{3.3}$$

Adachi et al. [1] have proved that $g^T y_2 = 0$ corresponds to the hard case, i.e., $\lambda_{opt} = -\alpha_n$ and $g \perp \mathcal{N}(A - \alpha_n B)$. Therefore, in the easy case, $g^T y_2 \neq 0$ is guaranteed, and (3.2) holds; furthermore, they have shown that $g^T y_2 = 0$ implies $y_1 = 0$. Numerically, in order to distinguish the hard and easy cases, they design a threshold $\tau < 1$ so that the case $\|y_1\| < \tau$ is recognized as the hard case, where $\tau$ is chosen as the level of square root of the machine precision.

Algorithm 1 is a modified form of Algorithm 2 in [1] that forms a framework of any eigenvalue-based algorithm for computing the solution $s_{opt}$ to (1.1) in the easy case.

---
**Algorithm 1** The eigenvalue-based algorithm in [1] (GEP)
---
1. (Consider the case $\lambda_{opt} = 0$) Solve $As_0 = -g$ for $s_0$ by the CG algorithm, and $s_{opt} = s_0$ if $\|s_0\|_B < \Delta$. Otherwise, go to step 2.
2. Use an eigensolver to compute the rightmost eigenvalue $\mu_1$ of the matrix pair $(M, \widetilde{B})$ and the corresponding eigenvector $y = (y_1^T, y_2^T)^T$.
3. If $\|y_1\| \leq \tau$ for a user-prescribed tolerance $\tau$, then treat TRS as hard case and solve it separately. Otherwise, $s_{opt} = -sign(g^T y_2) \Delta \frac{y_1}{\|y_1\|_B}$.
---

Adachi et al. have exploited IRA, i.e., the Matlab function eigs, in step 2 of this algorithm and made numerical experiments on varying sized TRS's to compare the overall efficiency of IRA with that of the non-restarted GLTR. In the experiments, they set stopping criteria for the relative residual norms in both eigs and GLTR to be the level of machine precision. The experiments have indicated that GLTR is a little more efficient. Comparisons have also been carried out in [34] for eigs, GLTR and their restarted GLTR as well as some others, where the stopping criterion for eigs was set to be the default value, i.e., the level of machine precision, and the ones for $B^{-1}$-relative residual norms of GLTR and others applied to (1.1) were $10^{-10}$. However, as we have commented in the introduction, such a comparison is problematic and unreasonable since the residual norms obtained by GLTR and an eigenvalue-based algorithm are defined differently, and for given user-prescribed stopping tolerances such as those in [1, 34], two kinds of algorithms may compute converged approximate solutions of TRS (1.1) with greatly different accuracy, so that the efficiency comparison is unfair and may be misleading. Therefore, it is vital to seek intimate relationships between these two kinds of residual norms and fix a user-prescribed stopping criterion for one kind and then decide the correct one for the other kind. Only in this way can we compute the approximate solution with very comparable accuracy and make a reasonable efficiency comparison. We shall solve this problem in Section 4.

**3.2. The IRA algorithm in [26].** In step 2 of Algorithm 1, one commonly uses the IRA algorithm to carry out the task [1, 34]. In this subsection, we shall briefly introduce some properties of the underlying Arnoldi method and the based IRA algorithm that are suited for computing Ritz approximations to the rightmost eigenpairs of a large matrix (pair).

Since $B$ is symmetric positive definite, the matrix $\widetilde{B}$ in (3.1) is too. We consider the eigenproblem of $(M, \widetilde{B})$:

$$(3.4) \qquad Mx_i = \mu_i \widetilde{B} x_i,$$

where $\mu_i$, $i = 1, 2, \ldots, 2n$, are the eigenvalues labeled in descending order of their real parts:

$$Re(\mu_1) \geq Re(\mu_2) \geq \cdots \geq Re(\mu_{2n}),$$

and the $x_i$ are corresponding eigenvectors of $\widetilde{B}$-norm one. By Theorem 3.1, we are concerned with $(\mu_1, y) := (\mu_1, x_1)$, the rightmost eigenpair of $(M, \widetilde{B})$, when solving TRS (1.1).

Given an initial $v_1$ of $\widetilde{B}$-norm one, the $k$-step Arnoldi process can be written in matrix form:

$$(3.5) \qquad MV_k = \widetilde{B} V_k H_k + h_{k+1,k} \widetilde{B} v_{k+1} e_k^T = \widetilde{B} V_{k+1} \widetilde{H}_k,$$

where $V_k = (v_1, v_2, \ldots, v_k) \in \mathbb{R}^{n \times k}$ and $V_{k+1} = (V_k, v_{k+1})$ are $\widetilde{B}$-orthonormal, $H_k \in \mathbb{R}^{k \times k}$ is upper Hessenberg, and $\widetilde{H}_k$ is the $(k+1) \times k$ upper Hessenberg matrix which is the same as $H_k$ except for an additional row whose only nonzero entry is $h_{k+1,k}$ in position $(k+1, k)$.

Denote by $(\mu_i^{(k)}, z_i^{(k)})$, $i = 1, 2, \ldots, k$, the eigenpairs of $H_k$ with $\|z_i^{(k)}\| = 1$, which are labeled as

$$Re(\mu_1^{(k)}) \geq Re(\mu_2^{(k)}) \geq \cdots \geq Re(\mu_k^{(k)}).$$

Then the Arnoldi method uses some of $\mu_i^{(k)}$, $i = 1, 2, \ldots, k$, called the Ritz values of the matrix pair $(M, \widetilde{B})$ with respect to the $k$-dimensional Krylov subspace

$$(3.6) \qquad \mathcal{K}_k(v_1, \widetilde{B}^{-1} M) \doteq span\{v_1, \widetilde{B}^{-1} M v_1, (\widetilde{B}^{-1} M)^2 v_1, \ldots, (\widetilde{B}^{-1} M)^{k-1} v_1\},$$

to approximate some exterior $\mu_i$, e.g., the rightmost one(s), and the corresponding eigenvectors $x_i$ of $(M, \widetilde{B})$ are approximated by the Ritz vectors of $\widetilde{B}$-norm one:

$$(3.7) \qquad x_i^{(k)} = V_k z_i^{(k)}.$$

The $\widetilde{B}^{-1}$-residual norm of any Ritz pair $(\mu_i^{(k)}, x_i^{(k)})$ satisfies

$$(3.8) \qquad \|(M - \mu_i^{(k)} \widetilde{B}) x_i^{(k)}\|_{\widetilde{B}^{-1}} = h_{k+1,k} |e_k^T z_i^{(k)}|,$$

which is used as a stopping criterion without forming $x_i^{(k)}$ until the convergence occurs.

However, an unfavorable aspect of the Arnoldi method is that one cannot know in advance how many steps will be required before the eigenpairs of interest are well approximated by the Ritz pairs. In order to recover Ritz vectors, one is obliged either to store $V_k$ or to recompute them. Therefore, the possible excessive storage and expensive computational cost limit the efficiency and applicability of the method. To this end, restarting is generally necessary.

Sorensen [26] has proposed an elegant implicit restarting technique: Suppose that some $k$ specific eigenpairs are of interest and the $(k+p)$-step Arnoldi process has been performed and one has computed $H_{k+p}$. Sorensen's means does not restart the Arnoldi process from the scratch. Rather, it implicitly restarts the Arnoldi process and realizes a new $(k+p)$-step Arnoldi process with an updated initial vector. It saves $k$ MVs with the matrix $\widetilde{B}^{-1}M$ at much less cost of $p$ implicitly shifted $QR$ iterations on the small matrix $H_{k+p}$ and naturally generates a new $k$-step Arnoldi process by truncating the implicitly shifted QR form. Then starting with iteration $k+1$, the standard Arnoldi process is performed to generate a new $(k+p)$-step Arnoldi process on the updated $(k+p)$-dimensional Krylov subspace. The resulting algorithm is the IRA (implicitly restarted Arnoldi) algorithm. One of the crucial points for the overall success of IRA is proper selection of the $p$ shifts involved. As a thumb, the shifts should be the best possible approximations to some unwanted eigenvalues of $(M, \widetilde{B})$. Within the framework of the Arnoldi method, Sorensen has suggested to use the $p$ unwanted eigenvalues of $H_{k+p}$ as shifts, called exact shifts [26]. Since then exact shifts have been used in all public IRA-based softwares, e.g., the Matlab function eigs. Algorithm 2 sketches IRA for computing the $k$ rightmost eigenpairs of $(M, \widetilde{B})$, where the superscript $H$ denotes the conjugate transpose of a matrix. We comment that Step 6 is performed using implicitly shifted QR iterations in real arithmetic in practical computations even if some shifts are complex conjugates.

**3.3. The IRRA algorithm in [16].** All the Arnoldi based algorithms, such as IRA, may have convergence problems [14, 18]: The Ritz vectors obtained by them may converge erratically and even may not converge even if the subspace contains sufficiently accurate approximations to the desired eigenvectors, while the corresponding Ritz values converge. To correct this deficiency, Jia [15, 16] has proposed a refined Arnoldi method and developed an IRRA (implicitly restarted refined Arnoldi) algorithm. The refined Arnoldi method is mathematically different from the Arnoldi method for extracting approximate eigenvectors from the underlying Krylov subspace, and the new eigenvector approximations, called the refined Ritz vectors, can be ensured to converge unconditionally provided that the subspace is sufficiently accurate. Moreover, the refined Ritz vectors are more accurate than the Ritz vectors [17].

In this subsection, we extend the refined Arnoldi method in [15] and the IRRA algorithm in [16] for the standard eigenvalue problem to the generalized eigenvalue problem of the matrix pair $(M, \widetilde{B})$. When adapting the refined Arnoldi method [15] to the generalized one, rather than using Ritz vectors $x_i^{(k)}$ as approximate eigenvectors, for each $\mu_i^{(k)}$ of interest we seek $\widetilde{B}$-normalized $\widetilde{x}_i^{(k)} \in \mathcal{K}_k(v_1, \widetilde{B}^{-1}M)$ satisfying the optimal condition

$$(3.9) \quad \left\|\left(M - \mu_i^{(k)}\widetilde{B}\right)\widetilde{x}_i^{(k)}\right\|_{\widetilde{B}^{-1}} = \min_{x \in \mathcal{K}_k(v_1, \widetilde{B}^{-1}M), \|x\|_{\widetilde{B}}=1} \left\|\left(M - \mu_i^{(k)}\widetilde{B}\right)x\right\|_{\widetilde{B}^{-1}},$$

and use it to approximate the desired eigenvector $x_i$. The approximation $\tilde{x}_i^{(k)}$ is called the refined Ritz vector, or more generally, the refined approximate eigenvector, from the underlying subspace associated with $\mu_i^{(k)}$.

The following results are from [15] for $\widetilde{B} = I$, and can be straightforwardly extended to the refined Arnoldi method applied to $(M, \widetilde{B})$ by making use of the

**Algorithm 2** IRA

1. Choose an initial vector $v_1$ with $\|v_1\|_{\widetilde{B}^{-1}} = 1$, the number $k$ of desired rightmost eigenpairs of $(M, \widetilde{B})$, the number $p$ of the shifts, and a stopping tolerance $\epsilon > 0$.
2. Perform the $k$-step Arnoldi process for $(M, \widetilde{B})$ with the starting vector $v_1$:
$$MV_k = \widetilde{B}V_k H_k + r_k e_k^T.$$
3. Apply $p$ additional steps of the Arnoldi process to obtain the $(k+p)$-step Arnoldi decomposition:
$$MV_{k+p} = \widetilde{B}V_{k+p}H_{k+p} + \widetilde{B}r_{k+p}e_{k+p}^T;$$
4. Compute the eigenpairs of $H_{k+p}$, and select its $k$ eigenvalues with the largest real parts to approximate the desired eigenvalues $\mu_i$, $i = 1, 2, \ldots, k$ of $M$ and the remaining $p$ ones as shifts $\sigma_1, \sigma_2, \ldots, \sigma_p$;
5. $Q = I$;
6. for $j = 1, 2, \ldots, p$
       QR factorization: $H_{k+p} - \sigma_j I = Q_j R_j$;
       $H_{k+p} = Q_j^H H_{k+p} Q_j$;
       $Q = QQ_j$;
   end
7. $H_k = H_{k+p}(1:k, 1:k)$;
8. $V_k = V_{k+p}Q(:, 1:k)$;
9. $r_k = H_{k+p}(k+1, k)V_{k+p}Q(:, k+1) + Q(k+p, k)r_k$;
10. If the right-hand sides of (3.8) drop below $\epsilon$ for $i = 1, 2, \ldots, k$, stop, and take $(\mu_i^{(k)}, x_i^{(k)})$ in (3.7) as approximations to $(\mu_i, x_i)$; otherwise, go to Step 3.

Arnoldi process (3.5).

THEOREM 3.2 ([15]). *Let $\widetilde{z}_i^{(k)}$ be the right singular vector of $\widetilde{H}_k - \mu_i^{(k)}\tilde{I}$ associated with $\sigma_{\min}\left(\widetilde{H}_k - \mu_i^{(k)}\tilde{I}\right)$, the smallest singular value of $\widetilde{H}_k - \mu_i^{(k)}\tilde{I}$, where $\tilde{I}$ is the same as $I$ except for an additional zero row. Then*

$$\widetilde{x}_i^{(k)} = V_k \widetilde{z}_i^{(k)}, \tag{3.10}$$

$$\left\|\left(M - \mu_i^{(k)}I\right)\widetilde{x}_i^{(k)}\right\|_{\widetilde{B}^{-1}} = \sigma_{\min}\left(\widetilde{H}_k - \mu_i^{(k)}\tilde{I}\right). \tag{3.11}$$

This theorem indicates that one can reliably and efficiently compute the refined Ritz vectors $\widetilde{x}_i^{(k)}$ by computing the singular value decomposition of a small sized $\widetilde{H}_k - \mu_i^{(k)}\tilde{I}$. Particularly, (3.11) is used as the stopping criterion without explicitly forming $\widetilde{x}_i^{(k)}$ until the convergence occurs.

The IRRA algorithm [16] applies the implicit restarting technique by Sorensen [26] to the Arnoldi process straightforwardly. A key contribution in [16] consists in providing a selection of better shifts, called the refined shifts, based on the information available during the refined Arnoldi method. The refined shifts have been theoretically shown and numerically justified to be better than the exact shifts used in IRA in the sense that they generate a better updated $(k+p)$-dimensional Krylov subspace during each cycle. Algorithm 3 sketches IRRA that computes the $k$ rightmost eigenpairs of $(M, \widetilde{B})$. We point out that Steps 5–6 are managed to be run only in real arithmetic and Step 8 is performed by using implicitly shifted QR iterations in real arithmetic even if there are complex conjugate shifts [18].

---

**Algorithm 3** IRRA

1. Choose an initial vector $v_1$ with $\|v_1\|_{\widetilde{B}^{-1}} = 1$, the number $k$ of desired rightmost eigenpairs of $(M, \widetilde{B})$, the number $p$ of the shifts, and a stopping tolerance $\epsilon > 0$.
2. Perform the $k$-step Arnoldi process for $(M, \widetilde{B})$ with the starting vector $v_1$:
$$MV_k = \widetilde{B} V_k H_k + r_k e_k^T.$$
3. Apply $p$ additional steps of the Arnoldi process to obtain the $(k+p)$-step Arnoldi decomposition:
$$MV_{k+p} = \widetilde{B} V_{k+p} H_{k+p} + \widetilde{B} r_{k+p} e_{k+p}^T;$$
4. Compute the eigenvalues $\mu_i^{(k)}$ of $H_{k+p}$ and the vectors $\widetilde{z}_i^{(k+p)}$ in Theorem 3.2 corresponding to the $k$ eigenvalues $\mu_i^{(k)}$ with the largest real parts, and define the matrix $Z_k = (\widetilde{z}_1^{(k+p)}, \widetilde{z}_2^{(k+p)}, \ldots, \widetilde{z}_k^{(k+p)})$;
5. Compute the full $QR$ factorization: $Z_k = (U_k, \widehat{U}_k) \begin{pmatrix} R_k & \times \\ O & \widehat{R}_k \end{pmatrix}$;
6. Compute the eigenvalues of the matrix $\widehat{U}_k^H H_{k+p} \widehat{U}_k$ as shifts $\sigma_1, \sigma_2, \ldots, \sigma_p$;
7. $Q = I$;
8. for $j = 1, 2, \ldots, p$
   $QR$ factorization: $H_{k+p} - \sigma_j I = Q_j R_j$;
   $H_{k+p} = Q_j^H H_{k+p} Q_j$;
   $Q = Q Q_j$;
   end
9. $H_k = H_{k+p}(1:k, 1:k)$;
10. $V_k = V_{k+p} Q(:, 1:k)$;
11. $r_k = H_{k+p}(k+1, k) V_{k+p} Q(:, k+1) + Q(k+p, k) r_k$;
12. If the right-hand sides of (3.11) drop below $\epsilon$ for $i = 1, 2, \ldots, k$, stop, and take $(\mu_i^{(k)}, \widetilde{x}_i^{(k)})$ in (3.10) as approximations to $(\mu_i, x_i)$; otherwise, go to Step 3.

---

For large scale TRS's, in Step 2 of Algorithm 1, we can use IRA and IRRA to compute only one rightmost eigenpair $(\mu_1, x_1)$ of $(M, \widetilde{B})$. Define $y = x_1$. We then recover $s_{opt}$ by the formula (3.3). The resulting algorithms are abbreviated as TRS_IRA and TRS_IRRA, respectively.

**4. Relationships between the TRS residual norms and the eigenvalue-based ones.** For $B \neq I$, by a change of variables:
$$\widehat{A} \leftarrow B^{-\frac{1}{2}} A B^{-\frac{1}{2}}, \quad \widehat{g} \leftarrow B^{-\frac{1}{2}} g, \quad \text{and} \quad c \leftarrow B^{\frac{1}{2}} s,$$
TRS (1.1) can be equivalently transformed into the standard TRS with $B = I$:

(4.1) $$\min_{\|c\| \leq \Delta} \widehat{q}(c) \quad \text{with} \quad \widehat{q}(c) = \widehat{g}^T c + \frac{1}{2} c^T \widehat{A} c,$$

Let $c_{opt}$ be the solution to (4.1) and $s_{\text{gltr}}^{(k)}$ and $c_{\text{gltr}}^{(k)}$ be the GLTR iterates applied to the subspaces $\mathcal{K}_k(g, B^{-1} A)$ and $\mathcal{K}_k(\widehat{g}, \widehat{A})$, respectively. Then it is shown in [19, 33] that GLTR generates the identical Lagrangian multipliers $\lambda^{(k)}$ for (1.1) and (4.1) and

$$s_{\text{gltr}}^{(k)} = B^{-\frac{1}{2}} c_{\text{gltr}}^{(k)}, \tag{4.2}$$

$$\|s_{\text{gltr}}^{(k)} - s_{opt}\|_B = \|c_{\text{gltr}}^{(k)} - c_{opt}\|, \tag{4.3}$$

$$q(s_{\text{gltr}}^{(k)}) - q(s_{opt}) = \widehat{q}(c_{\text{gltr}}^{(k)}) - \widehat{q}(c_{opt}), \tag{4.4}$$

$$\|(A + \lambda^{(k)} B) s_{\text{gltr}}^{(k)} + g\|_{B^{-1}} = \|(\widehat{A} + \lambda^{(k)} I) c_{\text{gltr}}^{(k)} + \widehat{g}\|, \tag{4.5}$$

More generally, independently of a specific solver for TRS (1.1) directly, let $\hat{c}$ be an approximate solution of (4.1), and define $\hat{s} = B^{-\frac{1}{2}} \hat{c}$ to be the approximate solution of (1.1). Then for a given $\hat{\lambda} \geq 0$ it is straightforward to show that

$$\|(A + \hat{\lambda} B) \hat{s} + g\|_{B^{-1}} = \|(\widehat{A} + \hat{\lambda} I) \hat{c} + \widehat{g}\|. \tag{4.6}$$

Define
$$\widehat{M} = \widetilde{B}^{-\frac{1}{2}} M \widetilde{B}^{-\frac{1}{2}}.$$

Then the eigenvalues of $\widehat{M}$ are identical to the eigenvalues $\mu_i$, $i = 1, 2, \ldots, 2n$ of $(M, \widetilde{B})$, and its eigenvectors $w_i$ are related to the eigenvectors $x_i$ of the matrix pair by $x_i = \widetilde{B}^{-\frac{1}{2}} w_i$, from which it follows that $\|x_i\|_{\widetilde{B}} = 1$ if $\|w_i\| = 1$ and vice versa. Mathematically, the Arnoldi method applied to $\mathcal{K}_k(v_1, \widetilde{B}^{-1} M)$ and $\mathcal{K}_k(\widetilde{B}^{-\frac{1}{2}} v_1, \widehat{M})$ compute the same Ritz values $\mu^{(k)}$ and the Ritz vectors $x^{(k)}$ and $w^{(k)}$, respectively, where $x^{(k)} = \widetilde{B}^{-\frac{1}{2}} w^{(k)}$, and their residual norms satisfy

$$\|(M - \mu^{(k)} \widetilde{B}) x^{(k)}\|_{\widetilde{B}^{-1}} = \|(\widehat{M} - \mu^{(k)} I) w^{(k)}\|. \tag{4.7}$$

More generally, independently of any specific eigensolver for the eigenproblem of $(M, \widetilde{B})$, let $(\hat{\lambda}, \hat{w})$ be an approximation to the rightmost eigenpair $(\lambda_{opt}, w)$ of $\widehat{M}$, and define $\hat{y} = \widetilde{B}^{-\frac{1}{2}} \hat{w}$ to be the corresponding approximate eigenvector of $(M, \widetilde{B})$. Then

$$\|(M - \hat{\lambda} \widetilde{B}) \hat{y}\|_{\widetilde{B}^{-1}} = \|(\widehat{M} - \hat{\lambda} I) \hat{w}\|. \tag{4.8}$$

Based on (3.3), we recover the approximation $\hat{s}$ of the solution $s_{opt}$ from $\hat{y}$ as follows: Let $\hat{y} = (\hat{y}_1^T, \hat{y}_2^T)^T$. Then the approximation $\hat{s}$ of $s_{opt}$ is recovered by

$$\hat{s} = -sign(g^T \hat{y}_2) \Delta \frac{\hat{y}_1}{\|\hat{y}_1\|_B}. \tag{4.9}$$

In terms of (4.5), (4.6) and (4.7), (4.8), clearly it suffices to consider the standard TRS (1.1) with $B = I$ and its mathematically equivalent eigenproblem when seeking possible intimate relationships between the residual norms of approximate solutions of TRS and those of approximate eigenpairs of $M$. Based on them, we are able to design reliable and proper stopping tolerances for GLTR and TRS_IRA, TRS_IRRA in order that they deliver converged approximate solutions of TRS (1.1) with very comparable accuracy, so that we can make a reasonable and fair comparison of the overall efficiency of the three algorithms.

Notice that a common stopping criterion for any iterative algorithm that solves (1.1) requires that the TRS residual norm satisfy

$$\|(A + \hat{\lambda} I) \hat{s} + g\| \leq tol_1 \tag{4.10}$$

for a given tolerance $tol_1$, where $\hat{\lambda}$ and $\hat{s}$ are obtained approximations to $\lambda_{opt}$ and $s_{opt}$, respectively. Correspondingly, for a large $M$, a stopping criterion for any eigenvalue-based iterative projection algorithm is

$$\left\|(M - \hat{\lambda}I)\hat{y}\right\| \leq tol_2, \tag{4.11}$$

with $\|\hat{y}\| = 1$ and $tol_2$ a given tolerance. Here we exclude those numerically backward stable algorithms such as the QR algorithm for small to moderate sized matrix eigenvalue problems, which automatically compute the eigenpair approximations at the level of machine precision in the sense of backward errors, i.e., $tol_2 = \|M\|\mathcal{O}(\epsilon_{\text{mach}})$ with $\epsilon_{\text{mach}}$ being the machine precision. Whenever (4.11) is met, we terminate the algorithm used, and recover the converged solution of TRS (1.1) as follows: From Theorem 3.1 and (3.3), $\hat{\lambda}$ is an approximation of $\lambda_{opt}$, and the converged approximate solution

$$\hat{s} = -\frac{\Delta^2}{g^T\hat{y}_2}\hat{y}_1 = -sign(g^T\hat{y}_2)\Delta\frac{\hat{y}_1}{\|\hat{y}_1\|}. \tag{4.12}$$

As we have stressed previously, in order to fairly compare the efficiency of an iterative algorithm for solving (1.1) directly and an eigenvalue-based one, a key premise is that two kinds of algorithms must compute the converged solutions with almost the same accuracy. But this crucial issue has been ignored in the literature, e.g., [1, 34]. Since one commonly measures solution accuracy by the residual norm of an approximate solution, the comparison premise becomes the following: Given an approximation $(\hat{\lambda}, \hat{y})$, suppose that it satisfies criterion (4.11) and $\hat{s}$ is recovered from $\hat{y}$ by the formula (4.12). Then how small is the corresponding TRS residual norm $\|(A + \hat{\lambda}I)\hat{s} + g\|$, that is, how should $tol_1$ in (4.10) be chosen? The other side of the problem is: Suppose that $(\hat{\lambda}, \hat{s})$ satisfies (4.10). Then what should $tol_2$ be chosen? Note that an eigenvalue-based algorithm can recover $\hat{s}$ from $\hat{y}$ but an algorithm, e.g., GLTR, that solves (1.1) directly does not provide an approximation $\hat{y}$ of the eigenvector $y$ of $(M, \widetilde{B})$. Therefore, it is only possible to decide $tol_1$ from $tol_2$, but not the other way around. To this end, we naturally attempt to select an appropriate value $tol_1$ for a given $tol_2$. This goal can be achieved, provided that we can find close relationships between $\|(A + \hat{\lambda}I)\hat{s} + g\|$ and $\|(M - \hat{\lambda}I)\hat{y}\|$.

We can present the following basic result.

THEOREM 4.1. *With the previous notation, we have*

$$\|(A + \hat{\lambda}I)\hat{s} + g\| \leq \frac{\Delta}{\|\hat{y}_1\|}\|(M - \hat{\lambda}I)\hat{y}\|. \tag{4.13}$$

*Proof.* Making use of (4.12), we obtain

$$\left\|(M - \hat{\lambda}I)\begin{pmatrix}\hat{y}_1 \\ \hat{y}_2\end{pmatrix}\right\|^2 = \left(\frac{\|\hat{y}_1\|}{\Delta}\right)^2\|(A + \hat{\lambda}I)\hat{s} + g\|^2 + \left\|\hat{y}_1 - (A + \hat{\lambda}I)\hat{y}_2\right\|^2. \tag{4.14}$$

Therefore, we have

$$\left(\frac{\|\hat{y}_1\|}{\Delta}\right)^2\|(A + \hat{\lambda}I)\hat{s} + g\|^2 \leq \|(M - \hat{\lambda}I)\hat{y}\|^2,$$

meaning that

$$\|(A+\hat{\lambda}I)\hat{s}+g\| \leq \frac{\Delta}{\|\hat{y}_1\|}\|(M-\hat{\lambda}I)\hat{y}\|. \quad \square$$

REMARK 4.1. *For $B \neq I$, based on (4.5), (4.6) and (4.7), (4.8), relation (4.13) becomes*

$$\|(A+\hat{\lambda}B)\hat{s}+g\|_{B^{-1}} \leq \frac{\Delta}{\|\hat{y}_1\|_B}\|(M-\hat{\lambda}B)\hat{y}\|_{B^{-1}}. \tag{4.15}$$

REMARK 4.2. *This theorem indicates that, given the stopping tolerance $tol_2 := tol$ in (4.11) for an eigenvalue-based algorithm, the TRS stopping tolerance $tol_1$ in (4.10) for an algorithm that solves (1.1) directly should be chosen as*

$$tol_1 = \frac{\Delta}{\|\hat{y}_1\|}tol. \tag{4.16}$$

*Notice that the two terms of the right-hand side in (4.14) are the squared norms of upper and lower parts of the left-hand side residual. They should be comparable in size. Therefore, the bound (4.13) is very sharp and almost an equality within roughly a multiple of two, showing that $tol_1$ defined by (4.16) is very reliable. This suggests us to use such a choice in GLTR for a given stopping tolerance in IRA and IRRA.*

REMARK 4.3. *For $\hat{y}_1$ very small, TRS is numerically in the hard case. (4.13) and the above remark show that the residual norm $\|(A+\hat{\lambda}I)\hat{s}+g\|$ may not be small, meaning that $\hat{s}$ is not a meaningful approximation of $s_{opt}$. For example, if $\|\hat{y}_1\|$ is as small as $\|(M-\hat{\lambda}I)\hat{y}\|$, $\|(A+\hat{\lambda}I)\hat{s}+g\| = \mathcal{O}(\Delta)$. This implies that the eigenvalue-based algorithms for computing the rightmost eigenpair of $(M, \widetilde{B})$ is not capable of solving the TRS in the hard case, as expected.*

REMARK 4.4. *In finite precision arithmetic, stopping criteria for relative residual norms must be used for reliability and general purpose. We terminate an eigenvalue-based algorithm whenever*

$$\frac{\|(M-\hat{\lambda}I)\hat{y}\|}{\|M\|_1} \leq tol$$

*is met for a user-prescribed tolerance $tol$, where $\|\cdot\|_1$ is the 1-norm of a matrix. Such criterion is independent of the scaling of $M$ and allows $tol$ to be at the level of $\epsilon_{\mathrm{mach}}$ and the premature or too late termination cannot occur. Correspondingly, the stopping criterion for algorithms that solve TRS directly becomes*

$$\frac{\|(A+\hat{\lambda}I)\hat{s}+g\|}{\|M\|_1} \leq \frac{\Delta}{\|\hat{y}_1\|}tol. \tag{4.17}$$

*Whenever $\|\hat{y}_1\|$ is fairly small for a modest $\Delta$, the relative residual norm of the recovered $\hat{s}$ cannot attain the level of $\epsilon_{\mathrm{mach}}$. This may be a limitation of eigenvalue-based algorithms for TRS.*

**5. The convergence of the residual norms obtained by GLTR, the Arnoldi method and the refined Arnoldi method for solving (1.1) with $B = I$.** In this section, we investigate the residual convergence of non-restarted GLTR, the Arnoldi method and the refined Arnoldi method for solving TRS (1.1) with $B = I$. The results are straightforwardly applicable to the TRS with $B \neq I$ by exploiting (4.6), (4.8) and the fact that the eigenvalues of $\widehat{M}$ and the pair $(M, \widetilde{B})$ are the same.

**5.1. The residual convergence of GLTR.** The authors have studied the residual convergence of the GLTR method in [19], where the following result is established.

THEOREM 5.1. *Suppose* $\|s_{opt}\| = \|s_{\text{gltr}}^{(k)}\| = \Delta$, *and define*

$$\eta_{k1} = \frac{\|g\|^2(\alpha_1 + \lambda_{opt})^2}{\|g\|^2 + (\alpha_1 + \lambda_{opt})^2\Delta^2}, \qquad \eta_{k2} = \frac{2(\alpha_1 + \lambda_{opt})^2}{\|g\|^2 + (\alpha_1 + \lambda_{opt})^2\Delta^2},$$

*where $\alpha_1$ and $\alpha_n$ are the rightmost and leftmost eigenvalues of $A$, respectively. Then*

$$(5.1) \quad \|(A + \lambda^{(k)}I)s_{\text{gltr}}^{(k)} + g\| \leq \left(\frac{4\eta_{k1}\Delta^2}{\|g\|} + 8(\alpha_1 + \lambda_{opt})\eta_{k2}\Delta^3\right)\left(\frac{\sqrt{\kappa} - 1}{\sqrt{\kappa} + 1}\right)^{2(k+1)}$$
$$+ 4\sqrt{\kappa}\Delta(\alpha_1 + \lambda_{opt})\left(\frac{\sqrt{\kappa} - 1}{\sqrt{\kappa} + 1}\right)^{k+1},$$

*where* $\kappa = \frac{\alpha_1 + \lambda_{opt}}{\alpha_n + \lambda_{opt}}$.

REMARK 5.1. *The second term of the right-hand side in (5.1) dominates the bound soon as $k$ increases, and $\|(A+\lambda^{(k)}I)s_{\text{gltr}}^{(k)}+g\|$ decreases at least like $\left(\frac{\sqrt{\kappa}-1}{\sqrt{\kappa}+1}\right)^{k+1}$. The numerical experiments have confirmed that the bound is very sharp and can be almost attainable.*

**5.2. The residual convergence of the Arnoldi method for TRS.** Recall that $(\lambda_{opt}, y) := (\mu_1, x_1)$ is the rightmost eigenpair of $M$ and $y = (y_1^T, y_2^T)^T$. In the Arnoldi method, $\lambda^{(k)} = \mu_1^{(k)}$ is the rightmost Ritz value of $M$ with respect to $\mathcal{K}_k(v_1, M)$ that is used to approximate $\lambda_{opt} = \mu_1$, and $y^{(k)} := x_1^{(k)}$ is the corresponding Ritz vector of $M$ that approximates $y$. Let

$$y^{(k)} = \begin{pmatrix} y_1^{(k)} \\ y_2^{(k)} \end{pmatrix}.$$

Then in terms of (4.12), an approximation $s_A^{(k)}$ of $s_{opt}$ is obtained by

$$(5.2) \qquad s_A^{(k)} = -\frac{\Delta^2}{g^T y_2^{(k)}} y_1^{(k)} = -\text{sign}(g^T y_2^{(k)})\Delta \frac{y_1^{(k)}}{\|y_1^{(k)}\|}.$$

Let $\pi_k = V_k V_k^T$ be the orthogonal projector onto $\mathcal{K}_k(v_1, M)$. Then $\pi_k M \pi_k$ is the restriction of $M$ to $\mathcal{K}_k(v_1, M)$ and its matrix representation is $H_k$ defined in (3.5). A direct application of Theorem 3.8 in [14] to our context gives the following result on the convergence of $\mu_1^{(k)}$ to $\mu_1 = \lambda_{opt}$.

LEMMA 5.2. *Suppose that $\|s_{opt}\| = \|s_A^{(k)}\| = \Delta$, and let $\varepsilon_k = \sin\angle(y, \mathcal{K}_k(v_1, M))$ with $y$ the eigenvector of $M$ associated with $\mu_1$. Then for $\varepsilon_k$ small it holds that*

$$(5.3) \qquad |\mu_1 - \mu_1^{(k)}| \leq \kappa(\mu_1^{(k)})\gamma_k\varepsilon_k + \mathcal{O}(\varepsilon_k^2),$$

*where $\kappa(\mu_1^{(k)})$ is the spectral condition number of $\mu_1^{(k)}$ and $\gamma_k = \|\pi_k M(I - \pi_k)\|$.*[1]

---

[1] In Theorem 3.8 of [14], $\sin\angle(y, \mathcal{K}_k(v_1, M))$ in the right-hand side of (5.3) is $\tan\angle(y, \mathcal{K}_k(v_1, M))$, but it is obvious that the sine and tangent can be replaced each other in the right-hand side when $\sin\angle(y, \mathcal{K}_k(v_1, M))$ becomes small.

Recall that $(\mu_1^{(k)}, z_1^{(k)})$ is the rightmost eigenpair of $H_k$ and $y^{(k)} = V_k z_1^{(k)}$. Let the columns of $Z_\perp^{(k)}$ be an orthonormal basis of the orthogonal complement of $\text{span}\{z_1^{(k)}\}$ so that $(z_1^{(k)}, Z_\perp^{(k)})$ is orthogonal. Then we have

$$\begin{pmatrix} (z_1^{(k)})^T \\ (Z_\perp^{(k)})^T \end{pmatrix} M_k (z_1^{(k)}, Z_\perp^{(k)}) = \begin{pmatrix} \mu_1^{(k)} & f_k^T \\ 0 & C_k \end{pmatrix}, \tag{5.4}$$

where $f_k^T = (z_1^{(k)})^T M_k Z_\perp^{(k)}$ and $C_k = (Z_\perp^{(k)})^T M_k Z_\perp^{(k)}$. Note that the eigenvalues of $C_k$ are the Ritz values other than $\mu_1^{(k)}$ of $M$ with respect to $\mathcal{K}_k(v_1, M)$.

In our notation, it is direct from Theorem 3.2 in [18] to obtain the following result.

LEMMA 5.3 ([18]). *With the previous notation, assume that* $\text{sep}(\mu_1, C_k) > 0$. *Then*

$$\sin \angle(y^{(k)}, y) \leq \left(1 + \frac{\|M\|}{\sqrt{1 - \varepsilon_k^2} \text{sep}(\mu_1, C_k)}\right) \varepsilon_k, \tag{5.5}$$

*where* $\text{sep}(\mu_1, C_k) = \sigma_{\min}(C_k - \mu_1 I)$ *is the separation of* $\mu_1$ *and* $C_k$.

Now we can establish the main result.

THEOREM 5.4. *With the previous notation, assume that* $\|s_{opt}\| = \|s_A^{(k)}\| = \Delta$, *and* $\text{sep}(\mu_1, C_k) > 0$. *Then*

$$\|(A + \mu_1^{(k)} I) s_A^{(k)} + g\| \leq c_k \varepsilon_k + \mathcal{O}(\varepsilon_k^2), \tag{5.6}$$

*where*

$$c_k = \frac{\Delta}{\|y_1^{(k)}\|} \left( \|M - \mu_1^{(k)} I\| \left(1 + \frac{\|M\|}{\sqrt{1 - \varepsilon_k^2} \text{sep}(\mu_1, C_k)}\right) + |\alpha| \kappa(\mu_1^{(k)}) \gamma_k \right)$$

*with* $\alpha = y^T y^{(k)}$.

*Proof.* Note that $y$ and $y^{(k)}$ are real with $\|y\| = \|y^{(k)}\| = 1$. Since

$$\sin \angle(y^{(k)}, y) = \|(I - yy^T) y^{(k)}\| = \|y^{(k)} - (y^T y^{(k)}) y\| = \|y^{(k)} - \alpha y\|$$

and $My = \mu_1 y$, we obtain

$$\begin{aligned} \|(M - \mu_1^{(k)} I) y^{(k)}\| &= \|(M - \mu_1^{(k)} I)(y^{(k)} - \alpha y) + \alpha(\mu_1 - \mu_1^{(k)}) y\| \\ &\leq \|(M - \mu_1^{(k)} I)(y^{(k)} - \alpha y)\| + \|\alpha(\mu_1 - \mu_1^{(k)}) y\| \\ &\leq \|M - \mu_1^{(k)} I\| \|y^{(k)} - \alpha y\| + |\alpha| |\mu_1 - \mu_1^{(k)}| \\ &= \|M - \mu_1^{(k)} I\| \sin \angle(y^{(k)}, y) + |\alpha| |\mu_1 - \mu_1^{(k)}|. \end{aligned}$$

Substituting bound (5.3) for $|\mu_1 - \mu_1^{(k)}|$ and bound (5.5) for $\sin \angle(y^{(k)}, y)$ into the above relation, we obtain

$$\|(M - \mu_1^{(k)} I) y^{(k)}\| \leq \left( \|M - \mu_1^{(k)} I\| \left(1 + \frac{\|M\|}{\sqrt{1 - \varepsilon_k^2} \text{sep}(\mu_1, C_k)}\right) + |\alpha| \kappa(\mu_1^{(k)}) \gamma_k \right) \varepsilon_k + \mathcal{O}(\varepsilon_k^2). \blacksquare$$

From Theorem 4.1,

$$\|(A + \mu_1^{(k)} I) s_A^{(k)} + g\| \leq \left( \frac{\Delta}{\|y_1^{(k)}\|} \right) \|(M - \mu_1^{(k)} I) y^{(k)}\|.$$

Combining the above two relations proves (5.6). □

This theorem shows how fast the TRS residual norm of the Arnoldi method converges depends on the size of $c_k$ and $\varepsilon_k$. As we can see, $c_k$ involves a few factors, where $\text{sep}(\mu_1, C_k) \approx \text{sep}(\mu_1^{(k)}, C_k)$ and $|\alpha| \approx 1$ as $y^{(k)} \to y$, and $\|y_1^{(k)}\| = \|z_1^{(k)}\|$. We shall discuss $\varepsilon_k$ later on.

**5.3. The residual convergence of the refined Arnoldi method for TRS.** Recall that $\widetilde{y}^{(k)} := \widetilde{x}_1^{(k)}$ is the refined Ritz vector of $M$ from $\mathcal{K}_k(v_1, M)$ to approximate $y$. Regarding the residual norms of the Ritz pair $(\mu_1^{(k)}, y^{(k)})$ and the refined Ritz pair $(\mu_1^{(k)}, \widetilde{y}^{(k)})$, Theorem 4.1 of [17] proves that the strict inequality holds:

$$\|(M - \mu_1^{(k)} I)\widetilde{y}^{(k)}\| < \|(M - \mu_1^{(k)} I)y^{(k)}\|,$$

provided that $(\mu_1^{(k)}, y^{(k)}) \neq (\mu_1, y)$, the exact rightmost eigenpair of $M$. This means that, in terms of residual norms, $\widetilde{y}^{(k)}$ is generally more accurate than $y^{(k)}$ as an approximation to $y$.

Let

$$\widetilde{y}^{(k)} = \begin{pmatrix} \widetilde{y}_1^{(k)} \\ \widetilde{y}_2^{(k)} \end{pmatrix}.$$

Then the refined Arnoldi method computes the approximation $s_{\text{RA}}^{(k)}$ of $s_{opt}$ by the formula

$$(5.7) \qquad s_{\text{RA}}^{(k)} = -\frac{\Delta^2}{g^T \widetilde{y}_2^{(k)}} \widetilde{y}_1^{(k)} = -\text{sign}(g^T \widetilde{y}_2^{(k)}) \Delta \frac{\widetilde{y}_1^{(k)}}{\|\widetilde{y}_1^{(k)}\|}.$$

The following result from [18, Theorem 7.1] estimates the residual norm $\|(M - \mu_1^{(k)} I)\widetilde{y}^{(k)}\|$.

LEMMA 5.5. *With the previous notation, we have*

$$(5.8) \qquad \|(M - \mu_1^{(k)} I)\widetilde{y}^{(k)}\| \leq \frac{\|M - \mu_1^{(k)} I\|\varepsilon_k + |\mu_1 - \mu_1^{(k)}|}{\sqrt{1 - \varepsilon_k^2}},$$

where $\varepsilon_k = \sin \angle(y, \mathcal{K}_k(v_1, M))$.

Combining (4.13) with (5.8) yields the following result on the TRS residual norm immediately.

THEOREM 5.6. *With the previous notation, we have*

$$(5.9) \qquad \|(A + \mu_1^{(k)} I)s_{\text{RA}}^{(k)} + g\| \leq \frac{\left(\|M - \mu_1^{(k)} I\| + \kappa(\mu_1^{(k)})\gamma_k\right)\Delta}{\|\widetilde{y}_1^{(k)}\|\sqrt{1 - \varepsilon_k^2}} \varepsilon_k + \mathcal{O}(\varepsilon_k^2).$$

Theorem 5.1 indicates that $\|(A + \lambda^{(k)} I)s_{\text{gltr}}^{(k)} + g\|$ decreases at least as fast as $\left(\frac{\sqrt{\kappa}-1}{\sqrt{\kappa}+1}\right)^{k+1}$, that is to say, $\kappa = \frac{\alpha_1 + \lambda_{opt}}{\alpha_n + \lambda_{opt}}$ is the main factor affecting the convergence of the residual norm in the GLTR method. Theorem 5.4 and Theorem 5.6 show that convergence of the TRS residual norms by the Arnoldi method and the refined Arnoldi method strongly relies on $\varepsilon_k = \sin \angle(y, \mathcal{K}_k(v_1, M))$.

A number of estimates on $\varepsilon_k$ have been available in [14, 15] for a general matrix unsymmetric $M$, which is directly applicable to our specific context and thus omitted. Unfortunately, it is impossible to draw any definitive conclusion on which of $\left(\frac{\sqrt{\kappa}-1}{\sqrt{\kappa}+1}\right)^{k+1}$ and $\varepsilon_k$ decays faster. As is seen from the results in [14, 15], among others, how fast $\varepsilon_k$ decays critically depends on the eigenvalue distribution of $M$, and the better separated is $\mu_1$ from $\mu_2$, the faster $\varepsilon_k$ decays. Except the fact that the rightmost eigenvalue $\mu_1 = \lambda_{opt}$, nothing is known on the other eigenvalues of $M$ and the relationships between them and the eigenvalues of $A$. Theoretically, it is well possible that the GLTR method converges faster in some cases but the converse may be true in other cases for the Arnoldi method and the refined Arnoldi method. Our later numerical experiments will confirm this.

TABLE 1

*Test matrices.*

| name | n | nonzeros | sparsity | application |
| --- | --- | --- | --- | --- |
| bcsstk29 | 13,992 | 619,488 | 0.3164% | Structural Problem |
| bcsstk33 | 8,738 | 591,904 | 0.7752% | Structural Problem |
| pcrystk02 | 13,965 | 968,583 | 0.4967% | Duplicate Materials Problem |
| garon2 | 13,535 | 373,235 | 0.2037% | Computational Fluid Dynamics Problem |
| lhr11 | 10,964 | 231,806 | 0.1928% | Chemical Process Simulation Problem |
| lhr11c | 10,964 | 233,741 | 0.1944% | Chemical Process Simulation Problem |
| barth5 | 15,606 | 61,484 | 0.0253% | Duplicate Structural Problem |
| nemeth01 | 9,506 | 725,054 | 0.8024% | Theoretical/Quantum Chemistry Problem Sequence |
| nemeth17 | 9,506 | 629,620 | 0.6968% | Subsequent Theoretical/Quantum Chemistry Problem |
| pkustk02 | 10,800 | 810,000 | 0.6944% | Structural Problem |
| nd3k | 9,000 | 3,279,690 | 4.0490% | 2D/3D Problem |
| coupled | 11,341 | 97,193 | 0.0756% | Circuit Simulation Problem |
| Pres_Poisson | 14,822 | 715,804 | 0.3258% | Computational Fluid Dynamics Problem |
| stokes64s | 12,546 | 140,034 | 0.0890% | Computational Fluid Dynamics Problem |
| tuma2 | 12,992 | 49,365 | 0.0292% | 2D/3D Problem |
| opt1 | 15,449 | 1,930,655 | 0.8089% | Structural Problem |
| ramage02 | 16,830 | 2,866,352 | 1.012% | Computational Fluid Dynamics Problem |
| Si5H12 | 19,896 | 738,598 | 0.1866% | Theoretical/Quantum Chemistry Problem |
| net25 | 9,520 | 401,200 | 0.4427% | Optimization Problem |
| pf2177 | 9,728 | 725,144 | 0.7663% | Optimization Problem |
| c-41 | 9,769 | 101,635 | 0.1065% | Optimization Problem |
| c-44 | 10,728 | 85,000 | 0.0739% | Optimization Problem |
| PGPgiantcompo | 10,680 | 48,632 | 0.0426% | Undirected Multigraph |
| wing_nodal | 10,937 | 150,976 | 0.1262% | Undirected Graph |

**6. Numerical experiments.** We now report numerical experiments to compare the overall performance of the three algorithms: TRS_IRA, TRS_IRRA and GLTR.

All the experiments were performed on an Intel Core (TM) i7 with CPU 3.6GHz, and 8 GB memory under the Microsoft Windows 10 64-bit system with the machine precision $\epsilon_{\mathrm{mach}} = 2.22 \times 10^{-16}$. We used the Matlab 2020a function eigs of IRA [26] and Matlab code of IRRA [18], called reigs, and the routine galahad_gltr of the GALAHAD library [9], a Fortran 90 code of GLTR.

Specifically, the Matlab code of TRSGEP[2] is available on the internet, where the Matlab function eigs is used to solve the eigenvalue problem. In our numerical experiments, we have taken different subspace dimensions $m = k + p = 20, 30, 50, 70$ in eigs and reigs. We take the stopping tolerance $tol = 10^{-12}$ for the relative residual norms of approximate eigenpairs of $(M, \widetilde{B})$ and the maximum number of restarts to be 600 in eigs and reigs. We use (4.17) to set the stopping tolerance of GLTR. Due to the paper length, we only report the results for the case that the subspace dimension is 30 since we have observed similar phenomena for the dimension $m = 20, 50, 70$. It should be noted that $m - (k + 3)$ shifts are used in the codes eigs and reigs. In our context, $k = 1$, and each cycle of eigs and reigs costs $m - 4$ MVs.

In galahad_gltr, we have used the double precision arithmetic, set the maximum number of iterations to be the problem dimension $n$, and performed the Lanczos process without reorthogonalization.

The matrices $A$ in our test problems are $A = G + G^T$, where the matrices $G$ are taken from the Suite Sparse Matrix Collection [3] with nine examples for testing a Lanczos method for the extreme Lorentz eigenvalue problem and fifteen matrices used in other applications [34], such as the fast factorization pivoting methods for sparse symmetric indefinite systems [25]. Some characteristics of these matrices are listed in Table 1, where *nonzeros* is the number of nonzero entries, and the sparsity is $nonzeros/n^2$. We take the matrix $B = I$ and the symmetric positive definite tridiagonal matrix $B = \text{tridiag}(1, 3, 1)$, respectively, which are used in [1, 34]. In TRS (1.1), we take the vector $g$ to be a unit length vector generated by the Matlab function randn(n,1), and the trust region radii are $\Delta = 1$ and $\Delta = 100$, respectively. We compute the relative residual norm

$$\text{(6.1)} \qquad \text{Res} \doteq \frac{\|(A + \hat{\lambda}B)\hat{s} + g\|_{B^{-1}}}{\|g\|_{B^{-1}}},$$

where $(\hat{\lambda}, \hat{s})$ is a converged approximation obtained by each of GLTR, TRS_IRA and TRS_IRRA using our stopping criteria.

The overall efficiency of each algorithm is measured by MVs. We should point out that, making use of the structures of $M$ and $B$ defined by (3.1), a matrix-vector product in eigs and reigs is two MVs in GLTR. That is, taking MVs in GLTR as reference, each iteration of eigs and reigs from $k + 3 = 4$ upwards costs two MVs. In the later tables, two MVs used in eigs and reigs precisely means one iteration. We shall report

$$\text{ratio} \doteq \frac{\text{MVs}_{\text{IRA}} - \text{MVs}_{\text{IRRA}}}{\text{MVs}_{\text{IRA}}},$$

where $\text{MVs}_{\text{IRA}}$ and $\text{MVs}_{\text{IRRA}}$ denote the MVs used by IRA and IRRA, respectively. This quantity measures the percentage of the overall efficiency improvement of TRS_IRRA over TRS_IRA.

Tables 2–3 report the results on the test problems with $B = I$ and $\Delta = 1, 100$, respectively. As we see from MVs in Tables 2–3, TRS_IRA and TRS_IRRA are robust and efficient, and their overall efficiency is, on average, lower than but comparable to that of GLTR for many of the problems. The numbers of iterations, which equal MVs/2 for TRS_IRA and TRS_IRRA, used by the former two algorithms and GLTR are roughly comparable for most of the test problems. For the problems garon2,

---

[2]TRSGEP is available at http://www.opt.mist.i.u-tokyo.ac.jp/~nakatsukasa/codes/TRSgep.m.

TABLE 2
*Numerical results on sparse $A = G + G^T$ with $G$ in Table 1 and $B = I$, $\Delta = 1$.*

| Matrix | GLTR | | TRS_IRA | | TRS_IRRA | | |
|---|---|---|---|---|---|---|---|
| | MVs | Res | MVs | Res | MVs | Res | ratio |
| bcsstk29 | 371 | $9.3e-08$ | 1226 | $2.0e-08$ | 1112 | $1.0e-08$ | 9.3% |
| bcsstk33 | 237 | $3.9e-08$ | 528 | $1.3e-08$ | 456 | $9.6e-09$ | 13.6% |
| pcrystk02 | 245 | $8.3e-08$ | 560 | $1.3e-08$ | 500 | $3.0e-08$ | 10.7% |
| garon2 | 213 | $1.25e-10$ | 316 | $2.9e-11$ | 292 | $2.6e-11$ | 7.6% |
| lhr11 | 137 | $2.8e-09$ | 320 | $1.8e-10$ | 296 | $1.7e-10$ | 7.5% |
| lhr11c | 137 | $2.8e-09$ | 320 | $3.6e-10$ | 296 | $1.9e-10$ | 7.5% |
| barth5 | 189 | $5.9e-09$ | 320 | $1.8e-10$ | 296 | $3.7e-10$ | 7.5% |
| nemeth01 | 247 | $1.9e-10$ | 488 | $2.8e-11$ | 396 | $5.1e-11$ | 23.0% |
| nemeth17 | 47 | $3.1e-12$ | 86 | $1.7e-12$ | 80 | $3.5e-12$ | 7.0% |
| pkustk02 | 181 | $2.1e-06$ | 422 | $1.7e-06$ | 368 | $1.5e-06$ | 12.8% |
| nd3k | 193 | $5.1e-10$ | 442 | $1.1e-10$ | 424 | $8.9e-11$ | 4.1% |
| coupled | 345 | $2.0e-08$ | 794 | $4.1e-09$ | 664 | $4.4e-09$ | 16.4% |
| Pres_Poisson | 127 | $8.0e-11$ | 264 | $1.1e-11$ | 232 | $7.6e-11$ | 12.1% |
| stokes64s | 169 | $1.4e-10$ | 404 | $7.6e-12$ | 360 | $1.5e-11$ | 10.9% |
| tuma2 | 127 | $2.4e-11$ | 200 | $2.8e-12$ | 180 | $7.9e-12$ | 9.0% |
| opt1 | 167 | $7.7e-08$ | 296 | $5.4e-08$ | 296 | $2.8e-08$ | 0.0% |
| ramage02 | 175 | $8.6e-08$ | 734 | $1.6e-08$ | 664 | $1.4e-08$ | 9.5% |
| Si5H12 | 215 | $3.0e-10$ | 372 | $4.2e-11$ | 360 | $6.5e-11$ | 3.2% |
| net25 | 31 | $3.0e-11$ | 128 | $3.2e-11$ | 108 | $2.5e-11$ | 15.6% |
| pf2177 | 285 | $1.4e-08$ | 660 | $4.2e-09$ | 620 | $7.5e-09$ | 6.1% |
| c-41 | 363 | $2.5e-07$ | 1088 | $1.4e-07$ | 1000 | $4.9e-07$ | 8.1% |
| c-44 | 267 | $2.7e-07$ | 788 | $4.0e-08$ | 758 | $2.1e-08$ | 3.8% |
| PGPgiantcompo | 101 | $6.7e-09$ | 216 | $6.5e-10$ | 194 | $9.3e-10$ | 10.2% |
| wing_nodal | 299 | $8.5e-09$ | 584 | $4.3e-10$ | 560 | $1.3e-10$ | 4.1% |

nemeth01 and tuma2 with $\Delta = 1$, TRS_IRA and TRS_IRRA use considerably fewer iterations than GLTR does. These results have confirmed that there is no considerable winner between GLTR and the eigenvalue-based algorithms TRS_IRA and TRS_IRRA. Regarding TRS_IRA and TRS_IRRA, it is clear that TRS_IRRA is more and can be considerably more efficient than TRS_IRA. A careful calculation shows that the average improvements are 9.15% for $\Delta = 1$ and 10.88% for $\Delta = 100$, respectively. Furthermore, the recorded Res's have important implications. When it is roughly equal to our stopping tolerance $10^{-12}$ for TRS_IRA and TRS_IRRA, $\Delta/\|\hat{y}_1\| = \mathcal{O}(1)$, meaning that $\|\hat{y}_1\|$ is not small and TRS is indeed in the easy case. If Res is much less than $10^{-12}$, then $\|\hat{y}_1\|$ is quite small, meaning that the solution accuracy of TRS is impaired though TRS is still in the easy case.

TABLE 3
*Numerical results on sparse $A = G + G^T$ with $G$ in Table 1 and $B = I$, $\Delta = 100$.*

| Matrix | GLTR | | TRS_IRA | | TRS_IRRA | | |
|---|---|---|---|---|---|---|---|
| | MVs | Res | MVs | Res | MVs | Res | ratio |
| bcsstk29 | 329 | $9.3e-04$ | 1224 | $8.6e-05$ | 1204 | $3.2e-05$ | 1.6% |
| bcsstk33 | 217 | $9.8e-05$ | 468 | $1.2e-05$ | 376 | $1.9e-05$ | 19.7% |
| pcrystk02 | 185 | $1.0e-03$ | 560 | $1.3e-04$ | 500 | $2.9e-04$ | 10.7% |
| garon2 | 1687 | $1.9e-06$ | 7194 | $3.5e-07$ | 5340 | $1.5e-07$ | 25.8% |
| lhr11 | 125 | $2.4e-05$ | 320 | $3.6e-06$ | 296 | $1.3e-06$ | 7.5% |
| lhr11c | 125 | $2.4e-05$ | 320 | $5.9e-06$ | 296 | $3.6e-06$ | 7.5% |
| barth5 | 151 | $5.3e-05$ | 320 | $2.3e-06$ | 296 | $1.7e-06$ | 7.5% |
| nemeth01 | 781 | $5.4e-06$ | 7132 | $1.2e-06$ | 6156 | $2.3e-06$ | 13.7% |
| nemeth17 | 261 | $6.6e-09$ | 1080 | $2.6e-09$ | 652 | $1.4e-09$ | 39.6% |
| pkustk02 | 149 | $2.6e-02$ | 378 | $1.5e-02$ | 368 | $3.3e-02$ | 2.7% |
| nd3k | 1297 | $3.9e-06$ | 10668 | $1.6e-07$ | 8328 | $3.3e-07$ | 21.9% |
| coupled | 281 | $4.0e-04$ | 722 | $2.1e-05$ | 652 | $7.1e-05$ | 9.7% |
| Pres_Poisson | 911 | $5.7e-07$ | 3218 | $1.1e-07$ | 2918 | $8.4e-07$ | 9.3% |
| stokes64s | 1991 | $1.8e-04$ | 32144 | $3.8e-05$ | 31260 | $3.6e-05$ | 2.8% |
| tuma2 | 487 | $4.8e-07$ | 1204 | $3.2e-08$ | 1196 | $2.7e-08$ | 0.7% |
| opt1 | 133 | $1.7e-03$ | 296 | $4.9e-04$ | 268 | $9.5e-04$ | 9.5% |
| ramage02 | 137 | $9.8e-04$ | 628 | $3.9e-04$ | 582 | $1.6e-04$ | 7.3% |
| Si5H12 | 197 | $3.5e-06$ | 420 | $5.1e-07$ | 416 | $2.9e-07$ | 1.0% |
| net25 | 29 | $3.0e-06$ | 118 | $1.1e-06$ | 108 | $1.0e-06$ | 8.5% |
| pf2177 | 245 | $1.5e-04$ | 724 | $3.4e-05$ | 620 | $2.4e-05$ | 14.4% |
| c-41 | 297 | $5.0e-03$ | 1088 | $1.5e-03$ | 1000 | $4.2e-03$ | 8.1% |
| c-44 | 217 | $1.3e-03$ | 788 | $3.0e-04$ | 690 | $2.2e-04$ | 12.4% |
| PGPgiantcompo | 89 | $4.5e-05$ | 230 | $5.7e-06$ | 200 | $5.4e-06$ | 13.1% |
| wing_nodal | 247 | $9.7e-05$ | 584 | $8.7e-06$ | 548 | $5.7e-06$ | 6.2% |

Since GLTR is non-restarted, its memory requirement becomes harsh and even unacceptable if too many iterations are needed unless the Lanczos basis vectors are stored in second memory, so that they can be called to produce the final solution. As a result, considerable communication time is needed whenever GLTR produces the final solution by calling the previously saved basis vectors. In contrast, TRS_IRA and TRS_IRRA have low fast storage requirement for fairly small projection subspaces, and it is easy and simple to implement them.

For $B = \mathrm{tridiag}(1, 3, 1)$, we observe that the results in Tables 4–5 are similar to those in Tables 2–3 for $B = I$, and the previous comments apply here. Besides garon2, nemeth01 and tuma2, the algorithms TRS_IRA and TRS_IRRA with $\Delta = 1$ also use much fewer iterations than GLTR for wing_nodal. In the meantime, a careful calculation shows that the average improvements are 8.73% for $\Delta = 1$ and 8.90% for $\Delta = 100$, respectively. A remarkable difference from Tables 2–3 is that all the Res's in Table 5 are much larger than $10^{-12}$, which implies that all the $\Delta/\|\hat{y}_1\|_B$ are quite large, or equivalently the $\|\hat{y}_1\|_B$ are small.

TABLE 4
Numerical results on sparse $A = G + G^T$ with $G$ in Table 1 and $B = \text{tridiag}(1, 3, 1)$, $\Delta = 1$.

| Matrix | GLTR | | TRS_IRA | | TRS_IRRA | | |
|---|---|---|---|---|---|---|---|
| | MVs | Res | MVs | Res | MVs | Res | ratio |
| bcsstk29 | 209 | $2.1e-09$ | 448 | $8.5e-10$ | 422 | $1.3e-09$ | 5.8% |
| bcsstk33 | 185 | $7.7e-09$ | 468 | $5.1e-09$ | 392 | $1.3e-09$ | 16.2% |
| pcrystk02 | 171 | $3.5e-09$ | 580 | $4.6e-09$ | 472 | $2.5e-09$ | 18.6% |
| garon2 | 255 | $4.3e-11$ | 424 | $1.8e-11$ | 372 | $2.4e-11$ | 12.3% |
| lhr11 | 147 | $1.3e-09$ | 392 | $5.1e-10$ | 344 | $1.7e-09$ | 7.5% |
| lhr11c | 147 | $1.3e-08$ | 422 | $9.4e-09$ | 368 | $4.8e-09$ | 12.8% |
| barth5 | 455 | $2.0e-08$ | 1080 | $3.7e-08$ | 1020 | $1.8e-08$ | 5.6% |
| nemeth01 | 199 | $1.9e-08$ | 360 | $5.7e-09$ | 320 | $7.8e-09$ | 11.1% |
| nemeth17 | 129 | $1.6e-11$ | 250 | $4.1e-12$ | 224 | $3.0e-11$ | 10.4% |
| pkustk02 | 187 | $2.0e-07$ | 378 | $1.2e-07$ | 368 | $1.4e-07$ | 2.7% |
| nd3k | 383 | $1.2e-10$ | 740 | $1.0e-10$ | 728 | $5.5e-11$ | 1.6% |
| coupled | 307 | $1.8e-06$ | 908 | $4.9e-06$ | 746 | $3.6e-06$ | 17.8% |
| Pres_Poisson | 143 | $1.9e-11$ | 298 | $4.2e-12$ | 268 | $1.3e-11$ | 10.1% |
| stokes64s | 243 | $8.0e-11$ | 512 | $9.2e-11$ | 476 | $7.2e-11$ | 7.0% |
| tuma2 | 299 | $4.3e-08$ | 630 | $3.2e-08$ | 580 | $3.2e-08$ | 7.9% |
| opt1 | 141 | $1.5e-08$ | 308 | $1.0e-09$ | 268 | $2.5e-09$ | 13.0% |
| ramage02 | 135 | $3.0e-08$ | 260 | $2.1e-08$ | 260 | $2.9e-08$ | 0.0% |
| Si5H12 | 455 | $4.5e-11$ | 854 | $3.4e-11$ | 782 | $4.2e-11$ | 8.4% |
| net25 | 63 | $3.8e-08$ | 116 | $7.5e-08$ | 108 | $6.5e-08$ | 6.9% |
| pf2177 | 115 | $6.9e-05$ | 260 | $2.1e-05$ | 236 | $3.1e-05$ | 9.2% |
| c-41 | 409 | $7.1e-07$ | 1932 | $5.2e-07$ | 1744 | $2.9e-07$ | 9.7% |
| c-44 | 421 | $1.7e-07$ | 1360 | $1.2e-07$ | 1310 | $1.8e-07$ | 3.7% |
| PGPgiantcompo | 105 | $6.7e-09$ | 238 | $2.5e-09$ | 228 | $5.2e-09$ | 4.2% |
| wing_nodal | 363 | $4.1e-08$ | 624 | $4.7e-09$ | 580 | $2.4e-09$ | 7.1% |

**7. Conclusion.** We have proposed and developed an eigenvalue-based formulation the refined Arnoldi method and IRRA algorithm to solve a large scale TRS's. The refined Arnoldi method can overcome the potential non-convergence deficiency of the Arnoldi method when computing eigenvectors. In order to compare any eigenvalue-based algorithm with the GLTR algorithm correctly and fairly, we have established the relationship between $\|(A+\hat{\lambda}B)\hat{s}+g\|_{B^{-1}}$ and $\|(M-\hat{\lambda}I)\hat{y}\|_{B^{-1}}$, where $\hat{s}$ is the approximation to $s_{opt}$ recovered from the converged eigenvector of the matrix pair $(M, B)$. This issue is important but has been ignored in the literature when comparing the efficiency and accuracy of an eigenvalue-based algorithm, such as IRA, and an iterative algorithm that solves TRS directly. Based on the obtained result, we can first preset a stopping tolerance for eigenvalue-based algorithms, e.g., IRA and IRRA, and then determine a reliable one for the algorithms that solve TRS directly, e.g., GLTR, so that the two kinds of algorithms can compute converged approximate solutions with very comparable accuracy. Therefore, we have provided a correct and fair comparison standard for the two kinds of algorithms. In the meantime, we have analyzed the residual norms obtained by GLTR and the Arnoldi method, the refined Arnoldi method for the equivalent eigenproblem of TRS, showing how they converge as the subspace dimension increases.

TABLE 5
*Numerical results on sparse $A = G + G^T$ with $G$ in Table 1 and $B = \text{tridiag}(1, 3, 1)$, $\Delta = 100$.*

| Matrix | GLTR | | TRS_IRA | | TRS_IRRA | | |
|---|---|---|---|---|---|---|---|
| | Mvs | Res | Mvs | Res | Mvs | Res | ratio |
| bcsstk29 | 175 | $2.3e-05$ | 448 | $2.9e-05$ | 422 | $1.9e-05$ | 5.8% |
| bcsstk33 | 165 | $5.8e-05$ | 422 | $1.5e-05$ | 372 | $4.3e-05$ | 11.8% |
| pcrystk02 | 139 | $8.2e-05$ | 632 | $1.9e-05$ | 472 | $1.3e-05$ | 25.3% |
| garon2 | 2121 | $4.4e-07$ | 12384 | $1.7e-07$ | 11188 | $9.4e-08$ | 9.7% |
| lhr11 | 133 | $1.3e-05$ | 370 | $5.2e-06$ | 360 | $2.0e-06$ | 2.7% |
| lhr11c | 133 | $1.3e-05$ | 422 | $7.3e-05$ | 380 | $3.3e-05$ | 10.0% |
| barth5 | 383 | $2.2e-04$ | 1008 | $4.3e-04$ | 928 | $1.2e-04$ | 7.9% |
| nemeth01 | 157 | $2.1e-04$ | 390 | $2.3e-05$ | 320 | $8.5e-05$ | 18.0% |
| nemeth17 | 157 | $7.0e-06$ | 384 | $4.8e-06$ | 368 | $3.1e-06$ | 4.2% |
| pkustk02 | 157 | $7.5e-04$ | 378 | $8.1e-04$ | 368 | $1.4e-04$ | 2.7% |
| nd3k | 2471 | $1.3e-06$ | 16544 | $9.5e-07$ | 16128 | $4.1e-08$ | 2.5% |
| coupled | 263 | $1.7e-02$ | 1230 | $1.7e-02$ | 828 | $1.8e-02$ | 32.7% |
| Pres_Poisson | 1025 | $2.2e-07$ | 4318 | $4.7e-08$ | 3648 | $2.1e-08$ | 15.5% |
| stokes64s | 111 | $3.1e-01$ | 9316 | $4.5e-01$ | 8640 | $4.5e-01$ | 7.3% |
| tuma2 | 223 | $4.3e-04$ | 630 | $4.1e-04$ | 584 | $4.7e-05$ | 7.3% |
| opt1 | 117 | $1.8e-04$ | 261 | $1.9e-04$ | 260 | $1.5e-04$ | 0.4% |
| ramage02 | 85 | $3.5e-04$ | 260 | $2.2e-04$ | 260 | $3.4e-04$ | 0.0% |
| Si5H12 | 483 | $4.4e-07$ | 854 | $3.8e-07$ | 782 | $1.7e-07$ | 8.4% |
| net25 | 51 | $8.7e-04$ | 116 | $9.2e-04$ | 108 | $7.7e-04$ | 6.9% |
| pf2177 | 55 | $2.0e-02$ | 260 | $7.8e-02$ | 236 | $3.5e-02$ | 9.2% |
| c-41 | 305 | $1.1e-02$ | 1744 | $7.5e-03$ | 1632 | $9.6e-03$ | 6.4% |
| c-44 | 331 | $2.5e-03$ | 1412 | $2.1e-03$ | 1308 | $4.7e-03$ | 7.4% |
| PGPgiantcompo | 89 | $5.5e-05$ | 256 | $1.1e-05$ | 232 | $6.8e-05$ | 9.4% |
| wing_nodal | 273 | $8.6e-04$ | 642 | $1.6e-04$ | 628 | $9.6e-05$ | 2.2% |

We have made extensive numerical experiments and confirmed that TRS_IRA and TRS_IRRA is comparable to GLTR. The results have also demonstrated that TRS_IRRA is more efficient than TRS_IRA and improves the latter by roughly 10% in terms of matrix-vector products.